\documentclass[12pt]{amsart}
\usepackage{amsfonts,amssymb}
\usepackage{mathrsfs}
\usepackage{array}
\usepackage[all,ps,cmtip,rotate]{xy} 
\usepackage[usenames,dvipsnames]{color}
\usepackage[hidelinks]{hyperref}
\usepackage{cleveref}
\usepackage{bm,mathtools}
\usepackage{enumitem}

\usepackage{cite}
\usepackage{amsmath}
\usepackage{mathrsfs}
\usepackage{geometry}
\usepackage{graphicx}
\usepackage{subfigure}
\usepackage{microtype}
\usepackage{geometry}

\newtheorem{define}{Definition}
\newtheorem{thm}[define]{Theorem}
\newtheorem{prop}[define]{Proposition}

\newtheorem{rem}[define]{Remark}
\newtheorem{cor}[define]{Corollary}
\newtheorem{lem}[define]{Lemma}

\def\CH{{\rm CH}}
\def\Griff{{\rm Griff}}
\def\alb{{\rm Alb}}
\def\pic{{\rm Pic}}

\def\AJ{{\rm AJ}}
\def\alg{{\rm alg}}
\def\hom{{\rm hom}}
\date{}
\title{Chow groups of Gushel-Mukai fivefolds}
\begin{document}

\author{Lin Zhou}
\address{\parbox{0.9\textwidth}{Beijing International Center For Mathematical Research, Peking University, 100871, Beijing, China}}
\email{{linzhoumath@pku.edu.cn}}
\begin{abstract}
\noindent We compute the Chow groups of smooth Gushel-Mukai varieties of dimension $5$.		
\end{abstract}
\maketitle

\section{Introduction}
Let $n\in \{3,4,5,6\}$. A Gushel-Mukai $n$-fold is a smooth Fano variety of dimension $n$, with Picard number $1$, index $n-2$ and degree $10$. By Mukai \cite{mukai1992curves} (together with the later result of Mella \cite{mella1997existence}), a smooth Gushel-Mukai (GM for short) $n$-fold $X$ can be obtained as a smooth dimensionally transverse intersection \cite[(1)]{debarre2020gusheli}
\begin{align}\label{equ-GM}
X=\rm {CGr}(2,V_5)\cap \mathbf{P}(W)\cap Q\subset \mathbf{P}(\mathbb{C}\oplus\bigwedge^2V_5),
\end{align}
where $V_5$ is a $5$-dimensional complex vector space, ${\rm CGr}(2,V_5)\subset\mathbf{P}(\mathbb{C}\oplus\bigwedge^2V_5)$ is the projective cone over the image of the Grassmannian variety ${\rm Gr}(2,V_5)$ under the  Pl\"ucker embedding, $W$ is a vector subspace of $\bigwedge ^2 V_5\oplus \mathbb{C}$ of dimension $n+5$, and $Q$ is a quadric hypersurface.

In this article, we are interested in the study of Chow groups of smooth GM fivefolds defined over the field of complex numbers $\mathbb{C}$.

The main result of this article is the following.
\begin{thm}\label{mainthm}
	Let $X$ be a smooth GM fivefold defined over $\mathbb{C}$. Then:
	\begin{enumerate}[label=(\roman*)]
		\item  The cycle class maps induce isomorphisms $\CH_i(X)\cong H_{2i}(X,\mathbb{Z}) \cong \mathbb{Z}$, for $i=0,1,4$, and $\CH_3(X)\cong H_{6}(X,\mathbb{Z}) \cong\mathbb{Z}\oplus \mathbb{Z}$.
		\item The $i$-th Griffiths group of $X$ vanishes for any $i\in \mathbb{N}$.
		\item The Abel-Jacobi map induces an isomorphism $\Phi:\CH_2(X)_{\alg}\to J^5(X)$, and we have the following splitting short exact sequence
		\begin{equation}\label{equ-exact}
		0\to J^5(X)\xrightarrow{\Phi^{-1}}\CH_2(X)\xrightarrow{\rm cl}\mathbb{Z}\oplus\mathbb{Z}\to 0,
		\end{equation}
		where $J^5(X)$ is the intermediate Jacobian of $X$, and $cl:{\CH}_2(X)\to \mathbb{Z}\oplus\mathbb{Z}$ is the cycle class map to $H^6(X,\mathbb{Z})\cong \mathbb{Z}\oplus\mathbb{Z}$ \cite[Proposition 2.6]{debarre2020gusheli}. 
	\end{enumerate}
\end{thm}
\begin{rem}
	\rm{The facts ${\CH_0(X)}\cong\mathbb{Z}$ and $ {\CH_4(X)}\cong\mathbb{Z}$ follow from the rationality of $X$ \cite[Proposition 4.2]{debarre2015gushel} and the fact that ${\pic(X)}=\mathbb{Z}H$ \cite[Lemma 2.29]{debarre2015gushel}, where $H$ is the hyperplane section class of $X\subset \mathbf{P}(W)$. Moreover,  Laterveer \cite[Proposition 3.1]{laterveer2020algebraic} proved that $\CH^j(X)\otimes_{\mathbb{Z}}\mathbb{Q}\to H^{2j}(X,\mathbb{Q})$ is an injection for $j\neq 3$, which led to the description of the $i$-th ($i\neq 2$) Chow group with rational coefficient.
			
	This article is inspired by the work of Fu and Tian \cite{fu20192} on cubic fivefolds. The main tools are decomposition of the diagonal which is initiated by Bloch-Srinivas \cite{bloch1974gersten}, and the relation between unramified cohomology and Chow groups (see \cite{colliot2012cohomologie}).
}
\end{rem}
\begin{center}
	\textbf{Acknowledge}
\end{center}
I would like to express my gratitude to my supervisors, Lie Fu and Zhiyu Tian, for their help. 
\section{Preliminaries}
\subsection{Gushel-Mukai fivefolds}
A GM fivefold $X$ is a smooth GM variety defined as in the introduction (\ref{equ-GM}) with $n=5$. There is a natural polarization on $X$, given by the restriction of the hyperplane class on $\mathbf{P}(W)$.

Thanks to the series of work of Debarre and Kuznetsov on GM varieties (\cite{debarre2015gushel},\cite{debarre2020gusheli},\cite{debarre2020gushelm},\cite{debarre2019gushel}), we know some of their properties, like their moduli spaces and period map (\emph{cf.} \cite{debarre2020gushelm}, \cite{debarre2019gushel})  and their Fano varieties of linear spaces (\emph{cf.} \cite{debarre2019gushel}). Here we list some facts about GM fivefolds that we are going to use.
\begin{prop}\label{prop-3}
	 Let $X$ be a smooth GM fivefold and $(V_6,V_5,A)$ the Lagrangian data associted with $X$ $(${cf.} \cite[Section 3]{debarre2015gushel}$)$. We set $Y_{A}^{\geq l}:=\{[v]\in \mathbf{P}(V_6)\mid{\rm dim}(A\cap(v\wedge\bigwedge^2V_6))\geq l\}$ and endow it with a natural scheme structure $(${cf.} \cite{o2006irreducible}$)$.
	 \begin{enumerate}[label=(\roman*)]
	 	\item \cite[Proposition 4.2]{debarre2015gushel} $X$ is rational;
	    \item \cite[Propositions 3.1 and 3.4]{debarre2019gushel} The group $H_*(X,\mathbb{Z})$ is torsion-free. The Hodge diamond of $X$ is as follows
	\begin{equation*}
	\begin{smallmatrix}
	&& &&&& 1 \\
	&& &&& 0 && 0  \\
	& &&&0&& 1 &&0 \\
	&& &0&& 0 && 0 &&0  \\
	&&0 &&0&& 2 &&0&&0 \\
	& 0&& 0 && 10 &&10&&0&&0  \\
	&&0 &&0&& 2 &&0&&0 \\
	&& &0&& 0 && 0 &&0  \\
	& &&&0&& 1 &&0 \\
	&& &&& 0 && 0  \\
	&& &&&& 1 
	\end{smallmatrix} 
	\end{equation*}
	\item \cite[Theorems 1.1 and 5.3]{debarre2020gusheli} Let $f_A:\widetilde Y^{\ge2}_{A}\to Y^{\ge2}_{A} $ be the double cover of $Y^{\ge2}_{A}$ branched along $Y^{\ge3}_{A}$ $($we call $\widetilde Y^{\ge2}_{A}$ the double EPW surface with respect to $A$ $(${cf.} \cite{debarre2019double}$)$$)$. There is an isomorphism of integral Hodge structures
	\begin{equation*}
	 H_1(\widetilde Y^{\ge2}_{A},\mathbb{Z})\cong H_5(X,\mathbb{Z})
	\end{equation*}
	which induces an isomorphism
	\begin{equation*}
	{\rm Alb}(\widetilde Y^{\ge2}_{A}) \cong {J}(X)
	\end{equation*}
	of principally polarized abelian varieties. In particular, if $Y_A^{\geq 3}=\varnothing$, this isomorphism can be induced by a subscheme $Z\subset X\times \widetilde{Y}_A^{\geq 2}$.
	\begin{equation}
	\label{eq:ajz-homology5}
	{\rm AJ}_{Z} \colon H_1(\widetilde Y^{\ge2}_{A},\mathbb{Z})
	{\stackrel{{}_{\scriptstyle\sim}}{\longrightarrow}} H_5(X,\mathbb{Z}).
	\end{equation}
\end{enumerate}
\end{prop}
\begin{cor}
	Let $X$ be a smooth GM fivefold, then the cycle class map induces an isomorphism $\CH_3(X)\cong H_6(X,\mathbb{Z})\cong \mathbb{Z}\oplus \mathbb{Z}$.
\end{cor}
\emph{Proof.} By \cite[Theorem 2.21]{voisin2019birational}, the Abel-Jacobi map $\CH_3(X)_{\hom}\to J^3(X)(\mathbb{C})$ is an isomorphism because of the rationality of $X$. But $H^3(X,\mathbb{Z})=0$ by \Cref{prop-3} (ii) and $X$ satisfies the integral Hodge conjecture of degree $4$ by \cite[Lemma 1.14]{voisin2019birational}, which implies that $\CH_3(X)\cong H^4(X,\mathbb{Z})\cong \mathbb{Z}\oplus \mathbb{Z}$.
\qed

We will use the followings results concerning the varieties of linear spaces contained in smooth GM fivefolds.
\begin{prop}\label{prop-4}
Let $X$ be a smooth GM fivefold and let $(V_6,V_5,A)$ be the Lagrangian data associted with $X$. We set $Y_{A,V_5}^{\geq l}:=\{[v]\in \mathbf{P}(V_5)\mid{\rm dim}(A\cap(v\wedge\bigwedge^2V_6))\geq l\}$ and endow it with a natural scheme structure $(${cf.} \cite{o2006irreducible}$)$. Let $Y_{A,V_5}:=Y_{A,V_5}^{\geq 1}$.
\begin{enumerate}[label=(\roman*)]
	\item \cite[Theorem 4.7]{debarre2019gushel} Let $F_1(X)$ be the Fano variety of lines of $X$. The map $\sigma\colon F_1(X) \to \mathbf{P}(V_5)$, which sends $[v\wedge V_3]$ to $[v]$ for $0\neq v\in V_3$, factors as
	\begin{equation}\label{equ-factor}
	F_1(X) \xrightarrow{\ \tilde\sigma\ } {\widetilde{\mathbf{P}(V_5)}} \xrightarrow{\ \ } \mathbf{P}(V_5), 
	\end{equation}
	where ${\widetilde{\mathbf{P}(V_5)}} \to \mathbf{P}(V_5)$ is the double cover branched along the sextic hypersurface $Y_{A,V_5} \subset \mathbf{P}(V_5)$ and $\tilde\sigma $ is a $\mathbf{P}^1$-bundle over the complement of  {the preimage of} $Y^{\ge 2}_{A,V_5} \cup \Sigma_1(X)$, where $\Sigma_1(X)$ is the \emph{kernel locus} of $X$ defined as \cite[Lemma 4.4]{debarre2015gushel}.
\item \cite[Lemma 4.3(b)]{debarre2019gushel}  Let $F_2(X)$ be the Fano variety of planes of $X$. Then $F_2(X)$ has two connected compontents, $F_2^{\sigma}(X)$ and $F_2^{\tau}(X)$, parametrizing the so-called $\sigma$-planes and $\tau$-planes. For $\sigma$-planes:\\
	$\bullet$ If $X$ is ordinary, or special with $\mathbf{p}_X \notin {\rm pr}_{Y,2}(E)$, the map $\sigma$ factors as
		\begin{equation*}
		F_2^\sigma(X) {\stackrel{{}_{\scriptstyle\sim}}{\longrightarrow}} \widetilde{Y}^{\ge 2}_{A,V_5} \longrightarrow Y^{\ge 2}_{A,V_5} \ensuremath{\lhook\joinrel\relbar\joinrel\to}\mathbf{P}(V_5), 
		\end{equation*}
		where $\widetilde{Y}^{\ge 2}_{A,V_5} \to Y^{\ge 2}_{A,V_5}$ is a double covering {of the curve $Y^{\ge 2}_{A,V_5}$} branched along $Y^{\ge 3}_{A,V_5}$.\\
		$\bullet$ If $X$ is special with $\mathbf{p}_X \in {\rm pr}_{Y,2}(E)$, the scheme $F_2^\sigma(X)$ is the union of {a double cover} $\widetilde{Y}^{\ge 2}_{A,V_5}$  and one  {double} or two  {reduced} components isomorphic to $\mathbf{P}^1$ and contracted by the map $\sigma$ onto $\Sigma_1(X)$.
	\item \cite[Lemma 5.8]{debarre2020gusheli} The Abel-Jacobi map 
		$$H_1(F^2_{\sigma}(X),\mathbb{Z})\to H_5(X,\mathbb{Z})$$ 
		is a surjective morphism of Hodge structures for general $X$.
\end{enumerate}
\end{prop}
\begin{cor}\label{cor-RCC}
	Let $X$ be a smooth GM fivefold. Then $F_1(X)$ is rationally chain connected. In particular, ${\CH_0}(F_1(X))=\mathbb{Z}$.
\end{cor}
\emph{Proof.} By \Cref{prop-4} (i), and Graber-Harris-Starr \cite[Corollary 1.3]{graber2003families}, it is enough to show that ${\widetilde{\mathbf{P}(V_5)}}$ is rationally chain connected. Suppose that $Y_{A,V_5} \subset \mathbf{P}(V_5)$ is defined by $f_0=0$, note that $Y_{A,V_5}$ is not a smooth sextic hypersurface, since $Y_{A,V_5}^{\geq 2}\subset {\rm Sing}(Y_{A,V_5})$ and $Y_{A}^{\geq 2}$ is an integral normal surface \cite[Theorem B.2]{debarre2015gushel}. Consider the flat family 
\begin{equation*}
\xymatrix{\widetilde{\mathbf{P}} \ar[d]^{p} \ar[r]&\mathbf{P}(1,1,1,1,1,3)\\
	\mathbb{H}&&	}
\end{equation*} 
where $\mathbb{H}:=\mathbf{P}(H^0(\mathbf{P}(V_5),\mathcal{O}(6)))$. The morphism $p$ is flat and for any $[f]\in \mathbb{H}$, $p^{-1}([f])$ is the hypersurface in $\mathbf{P}(1,1,1,1,1,3)$ defined by $x^2-f=0$ (here deg$(x)=3$). Therefore, the smoothness of $p^{-1}([f])$ is {equivalent} to the smoothness of the hypersurface $Z(f)\subset \mathbf{P}(V_5)$, in which case $p^{-1}([f])$ is a Fano manifold of index $2$. By Koll\'ar-Miyaoka-Mori \cite[Theorem 0.1]{kollar1992rational}, it is rationally chain connected (indeed is rationally connected) (\emph{cf.} \cite[Theorem 2.13, P.254]{kollar1999rational}). {By specialization},  $p^{-1}([f_0])={\widetilde{\mathbf{P}(V_5)}}$ is also rationally chain connected. 
\qed
\subsection{Unramified Cohomology}
Let $X$ be a smooth projective complex variety. Consider the identity continuous morphism 
$$\pi:X_{\rm cl}\longrightarrow X_{\rm Zar},$$
where $X_{\rm cl}$ is with the classical topology and $X_{\rm Zar}$ is with the Zariski topology. For any abelian group $A$ (e.g. $\mathbb{Z}$, $\mathbb{Q}$, $\mathbb{Q}/\mathbb{Z}$) and any positive integer $i$, we introduce a sheaf on $X_{\rm Zar}$ defined by the derived direct image of the constant sheaf $A$, that is, $\mathcal{H}^i(A):=R^i\pi_*A$. By definition, the \emph{$i-$th unramified cohomology} (\emph{cf.} \cite{bloch1974gersten},\cite{colliot2012cohomologie}) of $X$ with coefficients in $A$ is
$$H^i_{\rm nr}(X,A):=H^0(X_{\rm Zar},\mathcal{H}^i(A)).$$
It is clear that $H^i_{\rm nr}(X,A)=E_2^{0,i}$ in the Leray spectral sequence of ($\pi,A$).

On the other hand, the cohomology $H^{k}(X_{\rm cl},A)$ has a filtration by \emph{coniveau}: for any $c\in \mathbb{N}$,
$$N^cH^k(X,A)=\sum_{Z}^{} {\rm ker}(H^k(X,A)\to H^k(X\setminus Z,A))$$
where $Z$ runs through all the closed algebraic subsets of $X$ of codimension $\geq c$. From this filtration, we can get the coniveau spectral sequence (\emph{cf.} \cite[Section 3]{bloch1974gersten}):
$$'E_1^{p,q}=\bigoplus_{{Y}\in X^{(p)}}H^{q-p}(\mathbb{C}({Y}),A(-p))\Rightarrow N^{\bullet}H^{p+q}(X,A)$$
where $X^{(p)}$ is the set of integral subschemes of $X$ of codimension $p$.

Thanks to the work of Bloch–Ogus on the Gersten conjecture \cite[Proposition 6.4]{bloch1974gersten}, we have the following proposition due to Deligne:
\begin{prop}
	Let $X$ be a smooth projective complex variety. Starting from $'E_2$, the coniveau spectral sequence 
	$$'E_1^{p,q}=\bigoplus_{{Y}\in X^{(p)}}H^{q-p}(\mathbb{C}({Y}),A(-p))\Rightarrow N^{\bullet}H^{p+q}(X,A)$$
	coincides with the Leray spectral sequence associated with ($\pi,A$)
	$$E^{p,q}_2=H^p(X,\mathcal{H}^q(A))\Rightarrow H^{p+q}(X,A).$$
\end{prop}
Here we list some facts we will need about unramified cohomology and the Lerary spectral spectral of ($\pi,A$).
\begin{prop}\label{prop-nr}
	\begin{enumerate} [label=(\roman*)]
	\item  \cite[Proposition 3.3(iii)]{colliot2012cohomologie} If $X=\mathbf{P}^d_{\mathbb{C}}$, then
	\begin{equation*}
	H^p(X,\mathcal{H}^q(A))= 
	\begin{cases}
	0 & \text{if $p\neq q$},\\
	A & \text{if $p=q\leq d$}.
	\end{cases}
	\end{equation*}
	\item \cite[Corollary 7.4]{bloch1974gersten} There is a natural isomorphism 
	$${\CH }^p(X)/{\alg}\cong	H^p(X,\mathcal{H}^p(\mathbb{Z}))$$
	for any smooth projective complex variety $X$, where ${\CH} ^p(X)/{\alg}$ is the group of algebraic cycles of codimension $p$ modulo algebraic equivalence.
	\item \cite[Theorem 3.1]{colliot2012cohomologie} For $i,j\in \mathbb{N}$ and $n\in \mathbb{N}^*$, we have a short exact sequence
	\begin{equation}\label{eq-exact}
	0\to H^j(X,\mathcal{H}^i(\mathbb{Z}))\otimes\mathbb{Q}/\mathbb{Z}\to H^j(X,\mathcal{H}^i(\mathbb{Q}/\mathbb{Z}))\to H^{j+1}(X,\mathcal{H}^i(\mathbb{Z}))_{\rm tor}\to 0.
	\end{equation}
	\item \cite[Corollary 6.2]{bloch1974gersten} $H^p(X,\mathcal{H}^q(A))=0$ if $p>q$, for any smooth projective complex variety $X$.
\end{enumerate}
\end{prop}
As unramified cohomology is a birational invariant for smooth projective varieties \cite[Theorem 1.21]{voisin2019birational}. One deduce the following from \Cref{prop-3}(i).
\begin{cor}\label{cor-nr}
	For any smooth GM fivefold $X$ and any abelian group $A$, $H^0_{\rm nr}(X, A)\cong A$, and $H^i_{\rm nr}(X, A)=0$ for $i\neq 0$.
\end{cor}

\section{Chow group of $1$-cycles}
This section is devoted to the proof of ${\CH_1}(X)\cong \mathbb{Z}$. First, we are going to show that $\CH_1(X)\otimes_{\mathbb{Z}}\mathbb{Q}\cong\mathbb{Q}$, see Laterveer \cite{laterveer2020algebraic} for an alternative proof using the so-called Franchetta property. 
\begin{lem}\label{lem-connected by lines}
	Let $X$ be a smooth GM fivefold, for any two points $x,y\in X$, $x$ and $y$ can be connected by at most 4 lines.
\end{lem}
\emph{Proof.} 
Let $\gamma_X(x)=[U_2]$ and $\gamma_X(y)=[V_2]$ be two projective lines in $\mathbf{P}(V_5)$ associated with $x,y$, where $\gamma_X: X\to {\rm Gr}(2,V_5)\subset \mathbf{P}(\bigwedge^2V_5)$ is the Gushel map of $X$ defined as \cite[Section 2.1]{debarre2015gushel}.

\emph{Claim: $x$ and $y$ can be connected by lines if the lines $[U_2]$ and $[V_2]$ can be connected by lines in $\gamma_X(X)\subset {\rm Gr}(2,V_5)$.}\\
Proof of the claim: Let $F_1(X,x)\subset F_1(X)$ be the subscheme parametrizing lines passing through $x$, then $\sigma(F_1(X,x))=[U_2]$ by \cite[Proposition 4.1]{debarre2019gushel}. By the same reason, $\sigma(F_1(X,y))=[V_2]$. Note that $x_1$ and $x_2$ can be connected by at most two lines if $\gamma_X(x_1)\cap \gamma_X(x_2)\neq \varnothing$ (indeed if these two lines intersect at $[v]\in \mathbf{P}(V_5)$, $x_1$ and $x_2$ will be contained in a quadratic hypersurface $\rho_1^{-1}(v)$ in $\mathbf{P}^3$ or $\mathbf{P}^4$ by \cite[Proposition 4.5]{debarre2015gushel} ). Therefore, if there are lines $[U_2]=L_0, L_1,\dots,L_n=[V_2]$ in $\mathbf{P}(V_5)$ such that $L_i\in \gamma_X(X)$ and $L_i\cap L_{i+1}\neq \varnothing$ for $i=0,\dots,n-1$, $x$ and $y$ can be connected by $2n$ lines.

Therefore, it is enough to show that $[U_2]$ and $[V_2]$ can be connected by lines in $\gamma_X(X)\subset {\rm Gr}(2,V_5)$. Here dim$\gamma_X(X)=5$ as $\gamma_X$ is finite\footnote{deg$\gamma_X=1$ if $X$ is ordinary, deg$\gamma_X=2$ if $X$ is special in the sense of \cite[Section 2.5]{debarre2015gushel}.}. Note that fixed a point $[v]\in [U_2]$, there exist at least $2$-dimensional lines in $\gamma_X(X)$ passing $[v]$ because dim$\rho_1^{-1}(v)\geq 2$ and $\gamma_X$ is finite. Thus the universal family of lines in $\gamma_X(X)$ intersecting $[U_2]$ is at least $4$-dimensional, which must intersect with the projection line $[V_2]$ in $\mathbf{P}(V_5)$. In other words, there exist lines $[U_2]=L_0, L_1,L_2=[V_2]$ in $\gamma_X(X)$ such that $L_i\cap L_{i+1}\neq \varnothing$ for $i=0,1$.  
\qed
\begin{prop}\label{prop-Q rational}
	Let $X$ be a smooth GM fivefold, then ${\CH}_1(X)\otimes_{\mathbb{Z}}\mathbb{Q}\cong\mathbb{Q}$. Specifically, there exist an integer $N\neq 0$, such that for any $1$-cycle $C$, $NC=m[l]$ in ${\CH_1}(X)$ for some $m$, where $[l]$ is the class of a line in $X$.
\end{prop}
\emph{Proof.} By \Cref{cor-RCC}, all lines in $X$ have the same class in ${\CH}_1(X)$. Therefore, it suffices to show that for any curve $C$ in $X$, there is a positive integer $N$, such that $NC$ is rationally equivalent to sum of lines.

Let $g:F\to B$ parametrize the family of effective $1$-cycles with at most four lines components, and $u:F\to X$ be the evaluation map such that the map
$$m:F\times _B F\to X\times X$$
is generically finite of degree $N$ and surjective (this is possible, thanks to \Cref{lem-connected by lines}). By the arguement in \cite[Proposition 3.1]{tian2014one}, there are integers $m_i$'s and lines $F_i$'s in the fibers of $g$, such that $NC$ is rationally equivalent to $\sum_{i}^{}m_iu_*(F_i)$.
\qed
\begin{cor}\label{cor-decomposition of dia}
	Let $(X,H)$ be a polarized smooth GM fivefold. There exists a nonzero integer $m\in \mathbb{N}^*$ and a decomposition in $\CH_5(X\times X)$
	\begin{equation}\label{decomposition}
	m\triangle_X=m x\times X+m l\times h+Z',
	\end{equation}
	where $x$ is an arbitrary point of $X$, $l$ is an arbitrary projective line contained in $X$, $h:=c_1(H)$ is the hyperplane section class,  and $Z'$ is supported on $X\times T$ with $T$ a closed algebraic subset of pure dimension $3$ in $X$.
\end{cor}
\emph{Proof.} By the arguement of decomposition of the diagonal  \cite[Theorem 10.29]{voisin2003hodge}, there is a nonzero integer $m$ and a decomposition in $\CH_5(X\times X)$
\begin{equation}\label{decom}
m\triangle_X=Z_0+Z_1+Z',
\end{equation}
where  $Z_i$ is supported on $W_i \times W'_i$ with dim $W_i = i$ and dim $W'_i = n-i$ for $i=0,1$, and $Z'$ is supported on $X\times T$ with $T$ a closed algebraic subset of codimension $\geq 2$ in $X$. We suppose that $T$ is of pure dimension $3$. In particular, since $\CH_0(X)_{\mathbb{Q}}$ and $\CH_1(X)_{\mathbb{Q}}$ are generated by a point and a line in $X$ respectively, and ${\rm Pic}(X)=\mathbb{Z}H$ by \cite[Lemma 2.29]{debarre2015gushel}, we can write $Z_0=m_1x\times X$ and $Z_1=m_2 l\times h$ for some integers $m_1$ and $m_2$. 

It remains to show $m=m_1=m_2$.
Let (\ref{decom}) act by correspondences on the class of a point $x$ and the class of a line $l$. Note that $(x \times X)^*(x) =x$, $(x \times X)^*(l) = 0$, $(l \times h)^*(x) = 0$, $(l\times h)^*(l) = l$ and $Z'^*(x) = Z'^*(l) = 0$ by reason of dimension and $\triangle_X$ always acts as the identity, we get $mx = m_1 x$ and $ml = m_2 l$. Hence, $m = m_1 = m_2$.
\qed

Before proving ${\CH_1}(X)\cong \mathbb{Z}$, we recall the following well-known result.
\begin{lem}\label{lem-divisible}
	\cite[Lemma 7.10]{bloch1974gersten} The Chow group of algebraically trivial cycles is divisible.
\end{lem}
\begin{cor}\label{cor-griff}
	For a smooth GM fivefold $X$, let $\CH_1(X)_{\alg}$ and $\CH_1(X)_{\hom}$ be the subgroups of $\CH_1(X)$ of algebraically trivial cycles, resp. homologically trivial cycles. Then ${\Griff}_1(X)=0$, where ${\Griff}_i(X):=\CH_i(X)_{\hom}/\CH_i(X)_{\alg}$ is the $i$-th Griffiths group of $X$. In other words, $\CH_1(X)_{\alg}=\CH_1(X)_{\hom}$.
\end{cor}
\emph{Proof.}  By \cite[Proposition 1.30]{voisin2019birational}, the group ${\Griff}_1(X)$ is a birational invariant for smooth projective varieties. Note that $X$ is birationally isomorphic to $\mathbf{P}^5$ \cite[Proposition 4.2]{debarre2015gushel} and $\Griff_1(\mathbf{P}^5)=0$.
\qed
\begin{cor}\label{cor-alg equ}
	For any smooth GM fivefold $X$, $\CH_1(X)/\alg\cong \mathbb{Z}$, i.e. any curve $C$ in $X$ can algebraically equivalent to sum of lines. Here $\CH_1(X)/\alg$ is the group of $1$-cycles modulo algebric equvalence.
\end{cor}
\emph{Proof.} Indeed, $H^8(X,\mathbb{Z})=\mathbb{Z}l$, where $l$ is an arbitrary projective line in $X$. For any curve $C\subset X$, there is an integer $d$ such that $C-dl=0$ in $H^8(X,\mathbb{Z})$ (where $d$ is the degree of $C$). Therefore $(C-dl)\in \CH_1(X)_{\alg}$. That means $[C]=d[l]$ in $\CH_1(X)/\alg$.
\qed
\begin{prop}
	For any smooth GM fivefold $X$, ${\CH_1}(X)\cong \mathbb{Z}$.
\end{prop}
\emph{Proof.} By \Cref{cor-alg equ}, for any $1$-cycle $C$ of $X$, there is an integer $d$ such that $[C]-d[l]$ is an element in ${\CH_1}(X)_{\alg}$, where $l$ is an arbitrary line in $X$. But ${\CH_1}(X)_{\alg}$ is divisible, $[C]-d[l]=NC'$ in ${\CH_1}(X)$ for some $1$-cycle $C'$, where $N$ is the integer in \Cref{prop-Q rational}. Note that $NC'$ is rationally equivalent to, by \Cref{prop-Q rational}, a sum of lines $m[l]$. Thus $[C]=d[l]+NC'=d[l]+m[l]=(d+m)[l]$ in ${\CH_1}(X)$.
\qed

\section{The vaninshing of the Griffiths groups}

For a smooth GM fivefold $X$, we have shown that $\Griff_0(X)=\Griff_1(X)=\Griff_4(X)=0$ by \Cref{cor-griff} and the definition of Griffiths groups. We prove in this section that $\Griff_2(X)=\Griff_3(X)=0$. The main tool is unramified cohomology, which is a birational invariant for smooth projective varieties. 

We need the following lemma:
\begin{lem}
	\cite[Theorem 5.2 (2)]{debarre2019double} Suppose that $X$ is a smooth GM fivefold, and $A=A(X)$ is the Lagrangian subspace associted with $X$. There is a unique double cover $f_2:\widetilde{Y}_A^{\geq 2}\to {Y}_A^{\geq 2}$ branched along the finite set $Y_A^3$ {\rm (}note that $Y_A^3=Y_A^{\geq 3}$ as $A$ contains no decomposable vectors by \cite[Theorem 3.16]{debarre2015gushel} and \cite[Theorem B.2]{debarre2015gushel}{\rm)}, which is empty when $X$ is general. Moreover, the scheme $\widetilde{Y}_A^{\geq 2}$ is integral and normal, and it is smooth away from $f_2^{-1}(Y_A^3)$ and has ordinary double points along $f_2^{-1}(Y_A^3)$.
\end{lem}
\begin{prop}\label{prop-alge coniveau}
	Let $X$ be a smooth GM fivefold.
	\begin{enumerate}[label=(\roman*)]
		\item There exist a smooth projective curve $C$, and an algebraic $3$-cycle $W\subset C\times X$, such that the induced map $H^1(C,\mathbb{Z})\to H^5(X,\mathbb{Z})$ is a surjective morphism of Hodge structures. 
		\item $H^5(X,\mathbb{Z})$ is of algebraic coniveau $2$, i.e. $N^2H^5(X,\mathbb{Z})=H^5(X,\mathbb{Z})$.
	\end{enumerate}
	
\end{prop}
\emph{Proof.} (1) Consider the Hilbert closure $Z\subset X\times \widetilde{Y}_A^{\geq 2}$ defined in \cite[Lemma 5.2]{debarre2020gusheli}. By \cite[Theorem 5.3]{debarre2020gusheli}, if $Y_A^{\geq 3}=\varnothing$ (this condition holds for general $X$) so that $\widetilde{Y}_A^{\geq 2}$ is a smooth surface, the Abel-Jacobi map 
\begin{equation*}
{\AJ}_Z:H^3(\widetilde{Y}_A^{\geq 2},\mathbb{Z})\to H^5(X,\mathbb{Z})
\end{equation*}
is an isomorphism of integral Hodge structures. We can choose a very ample smooth divisor $i:C\to \widetilde{Y}_A^{\geq 2}$, then the Lefschetz hyperplane theorem (\emph{cf.} \cite[Theorem 1.23]{voisin2003hodge}) says that the composition
\begin{equation*}
H^1(C,\mathbb{Z})\to H^3(\widetilde{Y}_A^{\geq 2},\mathbb{Z})\to  H^5(X,\mathbb{Z})
\end{equation*}
is a surjection of Hodge structures. Indeed, we can just set $C:=\widetilde{Y}_{A,V_5}^{\geq 2}\cong F^2_{\sigma}(X)$ \cite[(26) and Lemma 5.8]{debarre2020gusheli}.

As for any $X$, consider the universal familys $\widetilde{\mathcal{Y}}_{\mathcal{A}}^{\geq 2}\to \mathbb{D}$, $\pi:\mathcal{X}\to \mathbb{D}$ and $\mathcal{Z}\to \mathbb{D}$ such that $\mathcal{Z}_t\subset \mathcal{X}_t\times \widetilde{\mathcal{Y}}_{\mathcal{A},t}^{\geq 2}$ is defined as \cite[Lemma 5.2]{debarre2020gusheli}. Here $\mathbb{D}$ is an analytical disk of small radius in $\mathbb{C}$ of center $0$ such that for all $t\in \mathbb{D}\setminus\{0\} $ the fiber $\mathcal{X}_t$ is general in the sense above and $\mathcal{X}_0\cong X$. In other words, we have the universal family:
\begin{equation}\label{equ-universal family}
\xymatrix
@M=6pt
{
	\mathcal{Z} \ar[d]^{p} \ar[r]^{q} & \mathcal{X} \ar[d]^{\pi} \\
	\widetilde{\mathcal{Y}}_{\mathcal{A}}^{\geq 2} \ar[r] & \mathbb{D}
}
\end{equation}

Choose a very ample smooth divisor $C\subset \widetilde{Y}_A^{\geq 2}$ away from the union of the singular points of $\widetilde{Y}_A^{\geq 2}$ and the finite subset of $\widetilde{Y}_A^{\geq 2}$ such that $Z\to \widetilde{Y}_A^{\geq 2}$ is not flat \cite[Lemma 5.2]{debarre2020gusheli}, and deform this smooth curve along $\mathbb{D}$, we get a flat family of smooth curves $\mathcal{C}\to \mathbb{D}$ (if has singular fiber, shrinking $\mathbb{D}$), the diagram (\ref{equ-universal family}) becomes the following when restricting to $\mathcal{C}$:
\begin{equation}
\xymatrix
@M=6pt
{
	\mathcal{W} \ar[d]^{p} \ar[r]^{q} & \mathcal{X} \ar[d]^{\pi} \\
	\mathcal{C} \ar[r] & \mathbb{D}
}
\end{equation}
Here $p$ and $\pi$ are flat. Note that $H^1(\mathcal{C}_0,\mathbb{Z})\cong H^1(\mathcal{C}_t,\mathbb{Z})$ for any $t\neq 0$, since $\mathcal{C}\to \mathbb{D}$ is an Ehresmann's fibration \cite[Theorem 9.3]{voisin2003hodge1}. And for any $t\neq 0$, the Abel-Jacobi map induced by $\mathcal{W}_t$, i.e. $H^1(\mathcal{C}_t,\mathbb{Z})\to H^5(\mathcal{X}_t,\mathbb{Z})$, is surjective. By continuity of local systems, the Abel-Jacobi map ${\AJ}_W:H^1(C,\mathbb{Z})\to H^5(X,\mathbb{Z})$ is a surjection, where $W:=\mathcal{W}_0$.

(2) Consider the algebraic cycle $W\subset C\times X$ constructed in (1). The Abel-Jacobi map $H_1(C,\mathbb{Z})\to H_5(X,\mathbb{Z})$ can factor as following:
\begin{equation*}
H_1(C,\mathbb{Z})\xrightarrow{q_*p^*} H_5(q(W),\mathbb{Z})\to H_5(X,\mathbb{Z})
\end{equation*}
Thus $H_5(q(W),\mathbb{Z})\to H_5(X,\mathbb{Z})$ is surjective. Equivalently, 
$H^5(X,\mathbb{Z})$ is supported on $q(W)$, which is an algebraic subset of codimension $k\geq 2$. Therefore, $H^5(X,\mathbb{Z})=N^kH^5(X,\mathbb{Z})\subset N^2H^5(X,\mathbb{Z})$.
\qed

\begin{prop}\label{prop-Griff}
	Let $X$ be a smooth projective complex GM fivefold.
	\begin{enumerate}[label=(\roman*)]
		\item ${\Griff}^2(X)=0$.
		\item ${\rm coker}\{H^5(X,\mathbb{Z}) \to H^1(X,\mathcal{H}^4(\mathbb{Z}))\}\cong {\Griff}^3(X)$.
	\end{enumerate}
	
\end{prop}
\emph{Proof.}
This proof is similar to the proof of \cite[Proposition 13]{fu20192}.
First, we have an exact sequence from the spectral sequence \cite[(8.2)]{bloch1974gersten}
\begin{equation}
H^3(X,\mathbb{Z})\to H^0(X,\mathcal{H}^3(\mathbb{Z}))\xrightarrow{d_2} H^2(X,\mathcal{H}^2(\mathbb{Z})) \xrightarrow{\rm cl} H^4(X,\mathbb{Z})
\end{equation}
where the first arrow is the composition $H^3(X,\mathbb{Z})\twoheadrightarrow E_{\infty}^{0,3}\hookrightarrow E_2^{0,3}= H^0(X,\mathcal{H}^3(\mathbb{Z}))$, the last arrow is the cycle class map. Therefore, by \Cref{cor-nr}, ${\ker}\{H^2(X,\mathcal{H}^2(\mathbb{Z})) \xrightarrow{\rm cl} H^4(X,\mathbb{Z})\}=0$. But this is exactly  $\Griff^2(X)$ by \Cref{prop-nr} (ii). (Indeed, the vanishing of $\Griff^2(X)$ can be obtained directly from \cite[Theorem 2.21]{voisin2019birational} as $\CH_0(X)\cong \mathbb{Z}.$).

Similarly, from the spectral sequence, we get an exact sequence \cite[(5)]{fu20192} (we use $H^4_{\rm nr}(X,\mathbb{Z})=0$ by \Cref{cor-nr}, and $H^3(X,\mathcal{H}^2(\mathbb{Z}))=0$ by \Cref{prop-nr}(iv).)
\begin{equation}\label{exact}
 H^2(X,\mathcal{H}^3(\mathbb{Z})) \to H^5(X,\mathbb{Z}) \to H^1(X,\mathcal{H}^4(\mathbb{Z}))\xrightarrow{d_2} H^3(X,\mathcal{H}^3(\mathbb{Z})) \to H^6(X,\mathbb{Z}).
\end{equation}
Hence ${\rm coker}\{H^5(X,\mathbb{Z}) \to H^1(X,\mathcal{H}^4(\mathbb{Z}))\}\cong {\rm ker} \{H^3(X,\mathcal{H}^3(\mathbb{Z})) \to H^6(X,\mathbb{Z})\}$. But the later is exactly  $\Griff^3(X)$ by \Cref{prop-nr} (ii).
\qed
\begin{cor}\label{cor-griff_23}
	Let $X$ be a smooth GM fivefold, $\Griff_2(X)$ is torsion free.
\end{cor}
\emph{Proof.} By \Cref{prop-nr} (iv) the first arrow of (\ref{exact}) is the composition of
 \begin{align*}
 \begin{split}
 H^2(X,\mathcal{H}^3(\mathbb{Z}))
 \to N^2H^5(X,\mathbb{Z})\to H^5(X,\mathbb{Z}).
 \end{split}
 \end{align*}
 It is surjective, since $X$ is of algebraic coniveau $2$ by \Cref{prop-alge coniveau}. Consequently, the second arrow in \Cref{exact} is zero, and
 $$H^1(X,\mathcal{H}^4(\mathbb{Z}))\cong\Griff^3(X)$$
by \Cref{prop-Griff} (ii).

 Let $i=4$ and $j=1$, the exact sequence (\ref{eq-exact}) implies $H^1(X,\mathcal{H}^4(\mathbb{Z}))_{\rm tor} \cong H^0(X,\mathcal{H}^4(\mathbb{Q}/\mathbb{Z}))$ as $H^0(X,\mathcal{H}^4(\mathbb{Z}))$ vanishes by \Cref{cor-nr}. But $H^0(X,\mathcal{H}^4(\mathbb{Q}/\mathbb{Z}))=0$ again by \Cref{cor-nr}. Therefore, ${\Griff_2(X)}_{\rm tor}\cong H^1(X,\mathcal{H}^4(\mathbb{Z}))_{\rm tor}=0$.
 \qed
\begin{prop}\label{prop-griff_2}
	Let $X$ be a smooth GM fivefold, $\Griff_2(X)=0$.
\end{prop}
\emph{Proof.} The proof is the same as \cite[Proposition 18]{fu20192}. By \Cref{cor-griff_23}, it is enough to show that $\Griff_2(X)$ is of torsion. 

We need the decomposition of the diagonal (\ref{decomposition}) in the proof of \Cref{cor-decomposition of dia}. Choose a desingularization $\widetilde{T}\to T$ equipped with a cycle $\widetilde{Z}'\in \CH^3(X\times \widetilde{T})$ whose image in $\CH^3(X\times T)$ is a multiple $m'Z'$ of $Z'$. We denote $\tilde{i}:\widetilde{T}\to X$ the composition of desingularization and inclusion. Note that $\widetilde{T}$ is not necessarily connected and we can write $\pic^0(\widetilde{T})=\prod \pic^0(\widetilde{T}_i)$ and $\alb(\widetilde{T})=\prod \alb(\widetilde{T}_i)$, where $\widetilde{T}=\coprod_i\widetilde{T}_i$ is the decomposition into connected components.

For any $2$-cycle $\alpha$, let the correspondences in (\ref{decomposition}) act on $\alpha$, then $\triangle_*(\alpha)=\alpha$, $(x\times X)_*(\alpha)=(l\times h)_*(\alpha)=0$ by the reason of dimension and  $m'Z'_*(\alpha)=\tilde{i}_*\circ\widetilde{Z}'_*(\alpha)$. Therefore, 
$mm'\alpha=\tilde{i}_*(\widetilde{Z}'_*(\alpha))$, i.e. $\tilde{i}_*\circ\widetilde{Z}'_*=mm'$ on $Z_2(X)$.

If $\alpha$ is homologically trivial, $\widetilde{Z}'_*(\alpha)$ is a $1$-cocycle in $\widetilde{Z}'$ which is homologically trivial. So $\tilde{i}_*(\widetilde{Z}'_*(\alpha))$ is algebraically trivial as ${\Griff^1}(\widetilde{T})=0$, so as $mm'\alpha$. Hence $\alpha$ is a torsion element in $\Griff_2(X)$.
\qed

\section{Proof of \Cref{mainthm} (iii)}
In the rest of this article, $X$ is a smooth GM fivefold defined on $\mathbb{C}$. In this section, we consider the Abel–Jacobi map
$$\Phi:\CH_2(X)_{\rm alg}\to J^5(X)(\mathbb{C})$$
and show that $\Phi$ is an isomorphism. Note that $\CH_2(X)_{\rm alg}$ is a divisible group by \Cref{lem-divisible}.
\begin{prop}\label{prop-isogeny}
	The Abel-Jacobi map $\Phi:\CH_2(X)_{\rm alg}\to J^5(X)$ is surjective with a finite kernel, i.e. is an isogeny.
\end{prop}
\emph{Proof.} 
Let us show the surjectivity of $\Phi$ first.
By \Cref{prop-alge coniveau} (i), we can consider the following commutative diagram by the compatibility between the Abel-Jacobi maps and the correspondence actions \cite[Theorem 12.17]{voisin2003hodge1}:
\begin{equation}
\xymatrix{{\CH_0(C)_{\rm alg}} \ar@{->}[r]^{\rm alb}_{\sim}\ar[d]^{W^*}&{ J}(C)\ar@{->>}[d]\\
\CH_2(X)_{\rm alg}\ar[r]^{\Phi}& J^5(X)	
}
\end{equation}
 the top horizontal arrow is surjective since $\CH_0(C)_{\rm alg}=\CH_0(C)_{\rm hom}$ and \cite[Lemma 12.11]{voisin2003hodge1}, and the right arrow is a surjection with connected fibers, because $H^1(C,\mathbb{Z})\to H^5(X,\mathbb{Z})$ is surjective. Therefore, the botton arrow $\Phi$ is a surjection.
 
 It remains to show the finitness of the kernel of $\Phi$. This part is the same as the proof of \cite[Proposition 18]{fu20192}. Let $\widetilde{T}$ and $\widetilde{Z}$ be as in the proof of \Cref{prop-griff_2}. We have a commutative diagram	
 \begin{equation}\label{equ-mm'}
 \xymatrix{
 	\CH^3(X)_{\rm alg}\ar[r]^{\Phi}\ar[d]^{\widetilde{Z'}_*}&J^5(X) \ar[d]^{[\widetilde{Z'}]_*}\\
 	\CH^1(\widetilde{T})_{\rm alg}\ar[r]^{\cong} \ar[d]^{\widetilde{i}_*} & \pic^0(\widetilde{T})\ar[d]^{[\widetilde{i}]_*}\\
 	\CH^3(X)_{\rm alg}\ar[r]^\Phi&J^5(X)
 }
 \end{equation}
 Similarly, (\ref{decom}) implies that the composites of the vertical arrows in (\ref{equ-mm'}) are both the multiplication by $mm'$. Now for all $\alpha\in {\ker}(\Phi)$, by the divisibility of $\CH^3(X)_{\rm alg}$, there exists a $\beta\in \CH^3(X)_{\rm alg}$, such that $\alpha = mm'\beta$. So $\Phi(\beta)\in J^5(X)[mm']$. Since the arrow in the middle of (7) is an isomorphism, we find that $\widetilde{Z'}_*(\beta)\in \CH^1(\widetilde{T})_{\rm alg}$$ [mm']$. From where,
 \begin{equation*}
 \alpha=mm'\beta=\widetilde{i_*}\circ\widetilde{Z'}_*(\beta)\in\widetilde{i_*}({\CH}^1(\widetilde{T})_{\alg} [mm'])
 \end{equation*}
 which is a finite set.
\qed 

\emph{Proof of \Cref{mainthm}(iii)}

This proof is similar to \cite[Proposition 19]{fu20192}. For any smooth $X$, there exists a commutative diagram for groups
\begin{equation*}
\xymatrix{{\CH_0(C)_{\rm alg}} \ar@{->}[r]^{\rm alb}_{\sim}\ar[d]^{W^*}&{J}(C)\ar@{->>}[d]^{}\\
	\CH_2(X)_{\rm alg}\ar[r]^{\Phi}& J^5(X)	
}
\end{equation*}
Let $A$ be the kernal of the right arrow, which is an abelian subvariety of $J(C)$. By \Cref{prop-isogeny}, as a morphism between abelian groups, $\Phi$ is surjective with a finite kernel. Since any morphism from a divisible group to a finite group is trivial, the image of $A$ in ${\ker} (\Phi)$ is trivial so that the image of the group  $\CH_0(C)_{\rm alg}$ in $\CH_2 (X)_{\rm  alg}$ is isomorphic to $J^5 (X)$. So we have a split of $\Phi$ and $\CH_2 (X)_{\rm  alg}\cong J^5(X)\oplus \ker (\Phi)$. However the group $\CH_2 (X)_{\rm  alg}$ is divisible (\Cref{lem-divisible}), and therefore $\ker (\Phi) = 0$, i.e. $\Phi$ is an isomorphism. Finally, the class map ${\CH}_2(X)\to H^6(X,\mathbb{Z})$ is surjective since $X$ satisfies the integral Hodge conjecture of degree 6 by \cite[Remark 4.2]{debarre2020gushel}. Moreover, we have the split exact sequence (\ref{equ-exact}).
\qed

\bibliographystyle{plain}
\bibliography{Chow_groups_of_Gushel-Mukai_5-folds}

\end{document}